\documentstyle[amscd]{amsart}
\numberwithin{equation}{section}
\theoremstyle{plain}
\newtheorem{thm}{Theorem}[section]
\newtheorem{cor}[thm]{Corollary}
\newtheorem{lem}[thm]{Lemma}
\newtheorem{prop}[thm]{Proposition}
\begin{document}
\title[Subshifts  from sofic shifts and Dyck shifts]{Subshifts from sofic shifts and Dyck shifts, \\
zeta functions and topological entropy }
\author{Kokoro Inoue}
\author{Wolfgang Krieger}
\begin{abstract}
We introduce a class of coded systems that we construct from  sofic systems and Dyck shifts and we study a class of subshifts that we obtain by excluding words of length two from Dyck shifts. We derive expressions for zeta functions  and topological entropy. We derive an expression for the zeta function of certain subshifts that we obtain by excluding words  from Dyck shifts and of certain subshifts that we obtain by excluding words  from the subshifts that are constructed from full shifts and Dyck shifts.
\end{abstract}
\maketitle
\def\det{{{\operatorname{det}}}}
\def\trace{{{\operatorname{trace}}}}
\def\card{{{\operatorname{card}}}}
   

\section{Introduction}
 Let $\Sigma$ be a finite alphabet. On $\Sigma^{\Bbb Z}$ there acts the shift that sends the point $(x_{i})_{i \in\Bbb Z } \in \Sigma^{\Bbb Z}$ into the point  $(x_{i+1})_{i \in\Bbb Z } \in \Sigma^{\Bbb Z}$. For a shift invariant set $Y\subset \Sigma^{\Bbb Z}$ we denote the set of periodic points in $Y$ by $P(Y)$ and we denote by $P_{n}(Y)$ the set of $p\in  P(Y) $ that have period $n\in \Bbb N$. The zeta function of a shift invariant set  $Y\subset \Sigma^{\Bbb Z}$ is defined by
$$
\zeta_{Y}(z) = e^{\Sigma_{n \in \Bbb N}\frac{1}{n}\card  (P_{n}(Y) )z^{n}}.
$$The dynamical systems that are given by the closed shift invariant subsets of $\Sigma^{\Bbb Z}$, with the restriction of the shift acting on them, are called subshifts. These are studied in symbolic dynamics. For an introduction to symbolic dynamics see \cite {Ki}  or  \cite {LM}. A word is called admissible for a subshift $X \subset  \Sigma^{\Bbb Z} $ if it appears somewhere in a point of $X$. We denote the language of admissible words of a subshift
$X \subset  \Sigma^{\Bbb Z} $ by ${\cal   L} (X)$, and we denote the set of words in ${\cal L} (X)$ of length $n $ by   ${\cal L}_n(X),n \in \Bbb N$. The $n$-block system of the subshift $X \subset  \Sigma^{\Bbb Z} $ is a subshift with alphabet $ {\cal L}_n(X)$, a topological conjugacy of the subshift onto its $n$-block system being given by the map
$$
x \to ((x_{j})_{i < j \leq i+n})_{i \in \Bbb Z}  \quad (x \in X).
$$
The topological entropy of a subshift $X \subset  \Sigma^{\Bbb Z} $ is given by
$$
h(X) = \lim_{n \to \infty} \frac{1}{n}\log \card  ( {\cal  L}_n(X)). 
$$

The edge shift of finite directed graph has as its language of admissible words the set of finite paths in the graph. With the adjacency matrix $A$ of the directed graph
the zeta function of the edge shift of the graph  is given by
$
\frac {1}{\det ({\bold 1} - Az)}
$
(see e.g. \cite [Theorem 6.4.6]   {LM}). A subshift  $X \subset  \Sigma^{\Bbb Z}$   is said to be of finite type, if there is a finite set ${\cal F}$ of words in the symbols of the alphabet $\Sigma$ such that $X$ is equal to the set of points in $ \Sigma^{\Bbb Z}$ in which no word in ${\cal F}$ appears. Sofic systems [W; LM, Chapter 3] are the images of subshifts of finite type under continuous shift commuting maps. The zeta function of sofic shifts is rational [LM, Theorem 6.4.8; BR, Section 3].
A finite directed labeled graph, in which every vertex has at least one incoming edge and at least one outgoing edge, presents a sofic system whose admissible words are the label sequences of finite paths in the graph. The labeling of a directed graph with label alphabet $\Sigma$ is called 1-right resolving if for every $\sigma \in \Sigma$ every vertex has at most one outgoing edge that carries the label $\sigma$. Every topologically transitive sofic system is canonically presented by a 1-right resolving irreducible finite directed graph that is known as its right Fischer automaton \cite {F}. 

For a fomal language ${\cal L}$ of words in the symbols of the alphabet $\Sigma$ we denote its 
generating function by $g({\cal L})$, and we denote by $Y({\cal L}  )$ the set of $x \in \Sigma^{\Bbb Z}$ such that there are indices $i_{k}, k \in \Bbb Z,$ such that $i_{k}<i_{k+1},k \in \Bbb Z,$ and 
\begin{align*}
 (x_{j})_{i_{k-1}\leq j < i_{k}} \in {\cal L},\qquad k \in \Bbb Z.  \tag {1.1}
\end{align*} 
The closure of $ Y( {\cal L}  )$ is a subshift that is called the coded system of $ {\cal L}$ \cite {BH}. A 
code 
in the symbols
 of the alphabet 
$\Sigma$ 
is said to be circular if for
$x \in P(Y({\cal C} ))$
the set 
$\{ i_{k},  k \in \Bbb Z \}$ of indices, such that (1.1) holds, is unique (\cite {BP}, section VII.1). For a circular code ${\cal C}$
$$
\zeta_{ Y( {\cal C } )   } = \frac {  1  }{ 1  - g( {\cal C}, z)},
$$
(see e.g. \cite [Proposition 4.7.11]{St}).

In this paper we construct coded systems that are close to previously described families of coded systems \cite {Kr, I, M}. We recall the relevant definitions and introduce notation.
Let $\Gamma$ be a finite set with more than one element. The set $\Gamma$  will remain fixed throughout the paper. We set
$$
N = \card (\Gamma). 
$$
We denote the generators of  the Dyck inverse monoid ${\cal D}_{N}$  (the polycyclic inverse monoid of \cite{NP}) by $ \gamma (-), \gamma (+)  , \gamma   \in  \Gamma $. These satisfy the relations
$$
\gamma(-)  \gamma^{\prime}(+) = \cases \bold 1, &\text {if $\gamma = \gamma^{\prime}$},\\
0, & \text{ if $ \gamma \neq \gamma^{\prime}, \qquad  \gamma , \gamma^{\prime}  \in  \Gamma $}.
 \endcases
$$
We recall the construction of the Dyck shifts \cite {Kr} and of the Motzkin shifts \cite{M, I}. 
The Dyck shift $D_{N}$ is the subshift with alphabet
 $ \{\gamma (-), \gamma (+)  :\gamma   \in  \Gamma \}$ 
and admissible words 
$(\gamma_{i})_{1 \leq i \leq I}, I \in \Bbb N,$ given by the condition 
\begin{equation*}
\prod_{1 \leq i \leq I} \gamma_{i} \neq  0. \tag {1.2}
\end{equation*}
A word $(\gamma_{i})_{1 \leq i \leq I} \in {\cal L} (  D_{N}), I \in \Bbb N,$ is called a Dyck word if 
\begin{equation*}
\prod_{1 \leq i \leq I} \gamma_{i} = \bold 1.  \tag {1.3}
\end{equation*}
The Motzkin shift  $M_{N }$ is the subshift with alphabet 
$\{\gamma (-), \gamma (+)  :\gamma   \in  \Gamma \}\cup \{ \bold 1  \} $
and admissible words $(\gamma_{i})_{1 \leq i \leq I}, I \in \Bbb N,$   also given by the condition (1.2), the Motzkin words being defined again by (1.3).
The Dyck shift $D_{ N}$ can be viewed as  the coded system of  the  Dyck code which is the set of  Dyck words that cannot be written as a non-trivial concatenation of Dyck words. The Motzkin shift $M_{ N}$ can be viewed as  the coded system of  the  Motzkin code which is defined as the set of  Motzkin words that cannot be written as a non-trivial concatenation of Motzkin words minus the set $\{\bold 1\}$.

In \cite {HI} a necessary and sufficient condition was given for an irreducible subshift of finite type to embed into a Dyck shift, and in \cite {HIK} this criterion was extended to a more general class of target shifts that were constructed by means of graph inverse semigroups. For the case that the graph inverse semigroup is the Dyck inverse monoid  ${\cal D}_{N}, N > 1,$  we recall the construction of these target shifts. Denote by ${\cal D}^{-}_{N} ({ \cal D}^{+}_{N} )$ the free semigroup with generators $\alpha_{-}(n)$ ($  \alpha_{+}(n)$), $ 1 \leq n \leq N$.  Let there be given a finite irreducible graph with vertex set ${\cal V}$ and edge set ${\cal E}$, together with a labeling map 
$$
\lambda: {\cal E} \to{ \cal D}^{-}_{N}  \cup   \{ { \bold 1 }\} \cup { \cal D}^{+}_{N} .
$$
The labeling map extends to paths 
$(e_{i})_{1 \leq i \leq I}$ in the graph  by
$$
\lambda ((e_{i})_{1 \leq i \leq I})) = \prod_{1 \leq i \leq I}\lambda (e_{i} ).
 $$ Assume that for all $V, V^{\prime} \in {\cal V}$ there exists a path $b$ in the graph $( {\cal V} ,{\cal E}    )$ that starts at $V$ and ends at $  V^{\prime} $ such that $\lambda (b) = \bold 1$. Also assume that for all 
 $V \in {\cal V}$ and for all $\gamma \in \Gamma$ there exist cycles $b^{(-)}$ and $b^{(+)}$ from $V$ to $V$ such that $\lambda(b^{(-)}) = \gamma(-)$ and $\lambda(b^{(+)}) = \gamma(+)$.
The labeled directed graph $( \cal V  ,{\cal E} ,   \lambda)$ presents a subshift 
$X({\cal V} , {\cal E} ,   \lambda)$ that is a subsystem of the edge shift of  the graph 
$({\cal V} ,{\cal E}    )$ with ${\cal L} (X({\cal V}  , {\cal E}  ,   \lambda) ) $ equal to the set of paths $b$ in the graph  $( {\cal V}  ,{\cal E}   )$  such that $  \lambda (b) \neq 0$. As in in  \cite {HIK} we call 
$X( {\cal V}  ,{\cal E}  ,   \lambda)$ a $ {\cal D}_{N} $-presentation. In this way the Dyck shift $D_{N}   , N > 1$, is presented by a graph with a single vertex and $2N$ loops that carry the labels $ \gamma(-),   \gamma(+), \gamma \in \Gamma$. The graph that presents the Motzkin shift has an additional loop that carries the label $\bold 1$. As in  \cite {HI} we say that a periodic point $p$ of $X{( \cal V}  ,{\cal E} ,   \lambda)$ has a negative (positive) multiplier if there exists an $i \in \Bbb Z$ such that, with a period $\pi$ of $p$, one has that $\lambda((x_{j})_{i \leq j < i+\pi})$ is in  ${\cal D}^{-}_{N}$($ { \cal D}^{+}_{N})$, and we say that $p$ is neutral if there exists an $i \in \Bbb Z$ such that $\lambda((x_{j})_{i \leq j < i+\pi}) = \bold 1$. One can attempt to construct a ${ \cal D}_{N}$-presentation from a given ${ \cal D}_{N}$-presentation 
 $X({\cal V}  , {\cal E}  ,   \lambda)$ by excluding for some $K > 1$ suitably chosen words of length K from $X( {\cal V}  ,{\cal E}  ,   \lambda)$. This amounts to excluding symbols from the $K$-block system of $X({\cal V} ,{\cal E}  ,   \lambda)$, which,  given the correct labeling, is itself a 
 ${\cal D}_{N}$-presentation. Provided one has not excluded too many words, one is left with a subshift that is still a ${ \cal D}_{N}$-presentation. The examples of ${ \cal D}_{N}$-presentations that were given in \cite {HIK} were obtained in this way, and the ${ \cal D}_{N}$-presentations that we will encounter in this paper arise in the same way.
  
In section 2 we describe a construction of subshifts that  generalizes  the constructions of the Dyck and Motzkin shifts. Also the Schr\"oder shifts appear here. The idea of the construction is to attach loops to a vertex of a finite irreducible 1-right resolving labeled graph, and to have these loops imitate the behavior of the generators of the Dyck inverse monoid. For instance, the labeled directed graph can be the Fischer automaton of an irreducible sofic system. In the case that the labeling of the directed graph is bijective this construction yields also  a ${\cal D}_{N}$-presentation of the constructed subshift.

After some preparations  in section 3, where we consider mappings that assign to the elements of a finite set non-empty subsets of the set, we study in section 4 the subshifts that are obtained by removing from the Dyck shifts $D_{N}$ words of the form $\gamma(+)\gamma^{\prime}(+)$, or, which is equivalent by symmetry, words of the form $\gamma(-)\gamma^{\prime}(-), \gamma \in \Gamma. $  

In sections 5 we consider two examples where we remove words of length three from  Dyck shifts  and Motzkin shifts. In section 6 we consider an example where we remove words of length three from the shifts that we construct in section 2 from the full shifts and the Dyck shifts. One checks that the 3-block system of the subshifts that are constructed in sections 5 and 6 are ${\cal D}_{N}$-presentations.

The aim is in all cases to obtain an expression for the zeta function of the subshifts and to determine their topological entropy. The zeta fuction we obtain by the same method as was used in \cite {Ke} to obtain the zeta function of the Dyck shifts  (see also \cite {I, KM}). In the subshift $X \subset \Sigma^{\Bbb Z}$ one identifies subshifts of finite type  or sofic shifts $Y^{-}, Y^{+}  \subset X$, and one identifies circular codes  ${\cal C}^{-}, {\cal C} ^{0}, {\cal C} ^{+} \subset {\cal L} (X)$  such that
 \begin{align*}
 P(X) = P( Y^{-} \cup Y^{+} ) \cup P( Y({\cal C}^{-} )\cup Y({\cal C}^{+} ) ),  \tag {1.4}
 \end{align*}
\begin{align*}
 P( Y^{-} \cup Y^{+} ) \cap P( Y({\cal C}^{-} )\cup Y({\cal C}^{+} ) ) = \emptyset, \tag {1.5}
 \end{align*}
\begin{align*}
 Y({\cal C}^{-} )\cap Y({\cal C}^{+} )=  Y({\cal C}^{0} ). \tag {1.6}
 \end{align*}
Then
 \begin{align*}
 \zeta_{X} = \dfrac { \zeta_{Y^{-}} \zeta_{Y({\cal C}^{-})}
  \zeta_{Y({\cal C} ^{+})} \zeta_{Y^{+}  }  } 
 { \zeta_{Y^{-} \cap Y^{+} }   \zeta_{Y({\cal C}^{0}     )  }   }. \tag {1.7}
\end{align*}
In case the subshift under consideration is a ${\cal D}_{N}$-presentation, $ P(Y({\cal C} ^{-}))   \cup  
  P(Y^{-})  $
coincides with the set of periodic poins of $X$ that are neutral or have a negative multiplier, and $  P(Y^{+})  \cup P(Y({\cal C}^{+}))  $ coincides with the set of periodic poins of $X$ that are neutral or have a  positive  multiplier.
 
The Sch\"utzenberger method \cite {Sc} (see also \cite {D}) has been for many years a standard method that has been routinely applied to solve enumeration problems as they arise in the computation of the generating functions of the circular codes 
from which we obtain zeta function and topological entropy of the coded systems that we consider here. In \cite {KM} the  Sch\"utzenberger method  was applied in the computation of zeta functions of a structurally significant class of coded systems, that includes the Dyck shifts. (Previously, in \cite {Kr} the topological entropy of the Dyck shift $D_{N}, N> 1,$ had been shown to be $\log (N+1)$  by methods from ergodic theory, and in \cite {Ke} the computation of the zeta function of the Dyck shifts was based on a probabilistic argument.) Zeta function and topological entropy of the Motzkin shifts were obtained in \cite {I}
 by the method of bijective correspondence. (See, however, the remark there on p. 3  on the  
 Sch\"utzenberger method.) In this paper we also apply the  Sch\"utzenberger method. However, in section  6, to obtain the generating function of the relevant code, we return to the method of bijective correspondence that was suggested in \cite{I}.

\section{A construction of subshifts}

Let there be given an irreducible finite directed labeled graph $G_{\circ}$ with a distinguished vertex  
$v _{\circ} $ and with labeling alphabet $\Sigma_{\circ}$ and labeling map 
$\lambda_{\circ}$. We assume that the labeling map 
$\lambda_{\circ}$ is 1-right resolving. It is allowed that  $G_{\circ}$ is the degenerate graph with the one vertex  $v _{\circ}$ and no edges. $\Sigma_{\circ}$ is allowed to contain only one symbol.  $G_{\circ}$  can present a topological Markov shift, or it can be 
the Fischer automaton of a topologically transitive sofic system. Let there also be given  $K_{\gamma}^{-}, K_{\gamma}^{+}  \in \Bbb N, \gamma \in  \Gamma $. We construct a labeled directed graph $G(\lambda_{\circ},v _{\circ}  ,(K_{\gamma}^{-} , K_{\gamma}^{+})_{\gamma \in \Gamma})$ from $G_{\circ}$ by attaching to the vertex $ v _{\circ}  $ directed loops that we name
$$
e(\gamma(-), k^{-}_{\gamma}) , \qquad 1 \leq  k^{-}_{\gamma} \leq K_{\gamma}^{-} ,
$$
$$
e(\gamma(+), k^{+}_{\gamma}) , \qquad 1 \leq  k^{+}_{\gamma} \leq K_{\gamma}^{+} ,
$$
and that we label by themselves, 
\begin{align}
&\lambda_{\circ} ( e(\gamma(-), k^{-}_{\gamma})) =  
 e(\gamma(-), k^{-}_{\gamma}) , \qquad 1 \leq  k^{-}_{\gamma} \leq K_{\gamma}^{-}  , \notag\\
&\lambda_{\circ} ( e(\gamma(+), k^{+}_{\gamma})) =  
 e(\gamma(+), k^{+}_{\gamma}) , \qquad 1 \leq  k^{+}_{\gamma} \leq K_{\gamma}^{+} , \qquad \gamma \in \Gamma . \notag
 \end{align}
Setting
$$
\varphi(\sigma_{\circ}) = {\bold 1},\  \   \  \qquad \  \sigma_{\circ} \in \Sigma_{\circ},
$$
and
$$
\varphi(e(\gamma(-), k^-)) = \gamma(-), \qquad 1\leq k^- \leq K^-,
$$
$$
\varphi(e(\gamma(+), k^+)) = \gamma(+), \qquad 1\leq k^+ \leq K^+,
$$ 
we call a path $e_{i}, 1 \leq i \leq I, I \in \Bbb N$, in  the graph $G(\lambda_{\circ},v _{\circ}  ,(K_{\gamma}^{-} , K_{\gamma}^{+})_{\gamma \in \Gamma})$ admissible if
$(\varphi( \lambda (e_{i}))_{1 \leq i \leq I}$ is an admissible word of the Motzkin shift.
We define a subshift $X(G(\lambda_{\circ},v _{\circ}  ,(K_{\gamma}^{-} , K_{\gamma}^{+})_{\gamma \in \Gamma}))$ with alphabet
$$
\Sigma_{\circ} \cup \{ e(\gamma(-), k^{-}_{\gamma}):  1 \leq  k^{-}_{\gamma} \leq K_{\gamma}^{-} \}  \cup 
 \{  e(\gamma(+), k^{+}_{\gamma}): 1 \leq  k^{+}_{\gamma} \leq K_{\gamma}^{+} \} 
$$
 by declaring a word as admissible for $X(G(\lambda_{\circ},v _{\circ}  ,(K_{\gamma}^{-} , K_{\gamma}^{+})_{\gamma \in \Gamma}))$ if it is the label sequence of an admissible path in the graph 
$G(\lambda_{\circ},v _{\circ}  ,(K_{\gamma}^{-} , K_{\gamma}^{+})_{\gamma \in \Gamma})$. We say that an admissible word of $X(G(\lambda_{\circ},v _{\circ}  ,(K_{\gamma}^{-} , K_{\gamma}^{+})_{\gamma \in \Gamma}))$ is literal-non-positive (literal-non-negative) if 
the images under $\varphi$ of all of it symbols are in  ${ \cal D}^{-}_{N}  \cup   \{  \bold 1 \} $ ($  \{  \bold 1 \} \cup  {\cal D}^{+}_{N} $). We say  that a set of admissible words of  $X(G(\lambda_{\circ},v _{\circ}  ,(K_{\gamma}^{-} , K_{\gamma}^{+})_{\gamma \in \Gamma}))$ is literal-uniform if all of its words are either literal-non-positive or literal-non-negative.
If the labeling $\lambda_{\circ}$ is bijective then 
$X(G(\lambda_{\circ},v _{\circ}  ,(K_{\gamma}^{-} , K_{\gamma}^{+})_{\gamma \in \Gamma}))$ is a ${\cal D} _N$-presentation of $X(G(\lambda_{\circ},v _{\circ}  ,(K_{\gamma}^{-} , K_{\gamma}^{+})_{\gamma \in \Gamma}))$, that uses on $G(\lambda_{\circ},v _{\circ}  ,(K_{\gamma}^{-} , K_{\gamma}^{+})_{\gamma \in \Gamma})$ the labeling $\lambda$ that labels the edges  of $G_{\circ} $ with $\bold 1$ and sets the label equal to the value of $\varphi$ for the other edges.
Let ${\cal R}$ be the set of paths in $G_{\circ}$ that start and end at the vertex $v_{\circ} $.
Let  ${\cal C}$ be the circular code  of paths $(e_{i})_{ 1 \leq i \leq I}, I \in \Bbb N,$ in the graph 
$G(\lambda_{\circ},v _{\circ}  ,(K_{\gamma}^{-} , K_{\gamma}^{+})_{\gamma \in \Gamma})$ that start and end at the vertex $v_{\circ} $ and that are such that the word $(\varphi(\lambda(e_{i}))_{ 1 \leq i \leq I} $  is in the Motzkin code.
Let  ${\cal C}^ -$ be the circular code that contains the paths in 
$G(\lambda_{\circ},v _{\circ}  ,(K_{\gamma}^{-} , K_{\gamma}^{+})_{\gamma \in \Gamma})$  that are concatenations of an admissible path with a literal-non-positive label sequence 
and a path in ${\cal C}$, let  ${\cal C}^0$ be the circular code that contains the paths in $G(\lambda_{\circ},v _{\circ}  ,(K_{\gamma}^{-} , K_{\gamma}^{+})_{\gamma \in \Gamma})$ that are concatenations of a path in ${\cal C}$ and of a path in ${\cal R}$, and let  ${\cal C}^ +$ be the circular code that contains the paths in $G(\lambda_{\circ},v _{\circ}  ,(K_{\gamma}^{-} , K_{\gamma}^{+})_{\gamma \in \Gamma})$ that are concatenations of a path in ${\cal C}$ and of an admissible path with a literal-non-negative label sequence.

We set
$$
K_- = \sum_{\gamma \in \Gamma}K^{-}_{\gamma},
 \quad K_+= \sum_ {\gamma \in \Gamma}K^{+}_{\gamma},
$$
$$
K = \sum_ {\gamma \in \Gamma}K^{-}_{\gamma}
K^{+}_{\gamma}.
$$

\begin{lem}\label{lem:2.1}
$$
g({\cal C}^{0}, z) =
\tfrac{1}{2}(1-\sqrt {1- 4Kg({\cal R}, z)^2z^2}).
$$
\end{lem}
\begin{pf}
It is
$$
g({\cal C}^{0}) = g({\cal C})g({\cal R}),
$$
and there is a counting argument (Sch\"utzenberger  method) that translates a set equation into the equation  
$$
\qquad\ \qquad \qquad \   \qquad\qquad \quad g({\cal C}, z) = 
\frac {Kg({\cal R}, z) z^2}
{1-g({\cal R}, z)
g({\cal C}, z)}. \qquad \qquad \qquad \ \ \qquad\qquad\qed
$$
\renewcommand{\qedsymbol}{}
\end{pf}

Denote by ${\cal D}^-$  the circular code of paths in $G(\lambda_{\circ},v _{\circ}  ,(K_{\gamma}^{-} , K_{\gamma}^+)_{\gamma \in \Gamma})$ that are obtained by letting one of the edges $e(\gamma(-), k^{-}_{\gamma}), 1 \leq  k^{-}_{\gamma} \leq K_{\gamma}^{-}$,  follow a path in ${\cal R}$ and denote by ${\cal D}^+$ be the circular code of paths in $G(\lambda_{\circ},v _{\circ}  ,(K_{\gamma}^{-} , K_{\gamma}^{+})_{\gamma \in \Gamma})$  that are obtained by letting  a path in ${\cal R}$  follow one of the edges $e(\gamma(+), k^{+}_{\gamma}),  1 \leq  k^{+}_{\gamma} \leq K_{\gamma}^{+}$. 

\begin{lem}\label{lem:2.2}
$$
g( {\cal C}^ -, z) = \frac{g({\cal C}^{0}, z) }{1- K_-g( {\cal R} , z)z},\quad g( {\cal C}^ +, z) = 
\frac{g({\cal C}^{0}, z) }{1- K_+ g( {\cal R} , z)z}.
$$
\end{lem}
\begin{pf}
Every path in ${\cal C}^ -$ can be written uniquely as a concatenation of a path that is a concatenation of paths in ${\cal D}^ -$ and of a path that is a concatenation of a path in ${\cal R}$ and a path in ${\cal C}$, and, symmetrically, every path in ${\cal C}^ +$ can be written uniquely as a concatenation of a path that is a concatenation of paths in ${\cal C}$ and a path in ${\cal R}$ and of a path that is a concatenation of paths in  ${\cal D}^ +$. Apply the Sch\"utzenberger method.
\end{pf}

Denote by ${\cal D}^-_{\circ}$  the circular code of words with a last symbol  $e(\gamma(-), k^{-}_{\gamma}), 1 \leq  k^{-}_{\gamma} \leq K_{\gamma}^{-}$,  that follows the label sequence of  a path in ${\cal R}$ and denote by ${\cal D}^+_{\circ}$ be the circular code of words
 with a first symbol  $e(\gamma(+), k^{+}_{\gamma}),  1 \leq  k^{+}_{\gamma} \leq K_{\gamma}^{+}$ that is followed by the label sequence of  a path in ${\cal R}$. 
Denote the sofic shift that is presented by $G_{\circ}$ by $Y^0$ and denote by  $ Y^-$($Y^+$) the coded system of ${\cal  D}^ -_{\circ} $($ { \cal D}^ + _{\circ}$). $ Y^-$ and $Y^+$  are sofic systems. Denote the adjacency matrix of $G_{\circ}$ by
 $A_{G_{\circ}}$. Also denote by $({\bold 1} - A_{G_{\circ}} z)^{ \langle  v_{\circ}\rangle}$ the matrix that is obtained by deleting in the matrix  $({\bold 1} - A_{G_{\circ}} z)$ the $v_{\circ}$-th row and the $ v_{\circ}$-th column.

\begin{thm}\label{thm:2.3}
\begin{align*}
g({\cal R}, z) = \frac{\det ({\bold 1} - A_{G_{\circ}} z)^{ \langle  v_{\circ}\rangle}    }{\det ({\bold 1} - A_{G_{\circ}} z)}, \tag {2.1}
\end{align*}
\begin{align*}
&\zeta_{X(G(\lambda_{\circ},v _{\circ}  ,(K_{\gamma}^{-} , K_{\gamma}^{+})_{\gamma \in \Gamma})}(z)= \\
&\frac{  2 \zeta_ {Y^0}(  1 +  \sqrt {1- 4Kg({\cal R}, z)^2z^2} ) }    { (\negthinspace 1- 2 K_{-}g({\cal R}, z)z +  
\sqrt {1- 4Kg({\cal R}, z)^2z^2} \negthinspace)\negthinspace(   \negthinspace 1- 2 K_{+}g({\cal R}, z)z +  
\sqrt {1- 4Kg({\cal R}, z)^2z^2} \negthinspace)    }.
\end{align*}
\end{thm}
\begin{pf}
(2.1) follows from an application of Cramer's rule.
$ Y^-,Y^+$, and ${\cal C}^{-},{\cal C}^{0},{\cal C}^{+}$, satisfy the relations (1.3-6). One has here 
$$
Y^-\cap Y^+=Y^0 .
$$
and also
$$
P( Y^-) = P({\cal D}^{-}_{\circ})\cup P(Y^0), \quad P( Y^+) = P(Y^0)\cup P({\cal D}^{+}_{\circ}). 
$$
One obtains the theorem by means of Lemma 2.1 and  Lemma 2.2 from (1.7) where one uses, that the 1-right resolving property of the labeling $\lambda_{\circ}$ implies that there is a one-to-one correspondence between the paths in ${\cal R}$ and their label sequences.
\end{pf}

\begin{cor}\label{cor: 2.4}
Under the assumption that $K_- \geq K_+$, the topological entropy of $X(G(\lambda_{\circ},v _{\circ}  ,(K_{\gamma}^{-} , K_{\gamma}^{+})_{\gamma \in \Gamma})$ is equal of the negative logarithm of the smallest positive root of
 $$
 zg({\cal R}, z) = \frac {K_-}{K_-^2+K}.
 $$
\end{cor}
\begin{pf}
One checks that
\begin{align*}
\limsup_{n \to \infty}   \frac{1}{n}\log \card  \ & P_n(X(G(\lambda_{\circ},v _{\circ}  ,(K_{\gamma}^{-} , K_{\gamma}^{+})_{\gamma \in \Gamma})) = \\
&  \lim_{n \to \infty} \frac{1}{n}\log \card  \ {\cal  L}_n(X(G(\lambda_{\circ},v _{\circ}  ,(K_{\gamma}^{-} , K_{\gamma}^{+})_{\gamma \in \Gamma})),
 \end{align*}
and one  applies Theorem 2.1 and Theorem 2.3.
\end{pf}

We see the Dyck shift $D_N$ here as the special case  $Y = \emptyset,g({\cal R})= 1 $,  $K^-_\gamma =  K^+_\gamma = 1, \gamma \in \Gamma$.
Also, as an example for $G_{\circ}$ take a bouquet $G_{\bigcirc}(J,Q)$  of $J$ circles, each circle made up of $Q$ edges, $J, Q \in \Bbb N$: Besides the vertex $v_{\bigcirc}, G_{\bigcirc}$ has vertices $v_{j, q}, 1\leq j \leq J,   1\leq q< Q $, and there is an edge from the vertex $v_{\bigcirc}$ to each of the vertices $v_{j, 1}, 1\leq j \leq J $, there is an edge from the vertex $v_{j, q}$ to the vertex $v_{j, q+1}, 1\leq j \leq J,   1\leq q < Q-1 $, and an edge from each of the vertices $v_{j, Q-1}, 1\leq j \leq J, $ to  $v_{\bigcirc}$. With $\lambda_{\bigcirc}(J,Q)$ the labeling map that labels the edges of $G_
{\bigcirc}$ by themselves we write $G_{\bigcirc}(N,J,Q) $ for
$G(\lambda_{\bigcirc}(J,Q), v_{\bigcirc}, (1,1)_{\gamma \in \Gamma})$. $G_{\bigcirc}(N,1,1) $ is the Motzkin shift $M_N$. The subshifts $X(G_{\bigcirc}(N,J,1))$ were introduced in \cite  {I}  where their zeta function  \cite [Proposition 2.3] {I} and topological entropy \cite [Corollary 2.2] {I} were determined.

\begin{cor}\label{cor: 2.5}
 $$
 \zeta_{X(G_{\bigcirc}(N,J,Q))}(z)=
 \frac{2(1-Jz^Q + \sqrt{(1-Jz^Q)^2 - 4Nz^2})}
 {(1 - Jz^Q - 2Nz + \sqrt{(1-Jz^Q)^2 - 4Nz^2})^2},
\qquad N > 1.
$$
\end{cor}
\begin{pf}
In this case 
$$
\qquad \qquad \qquad\qquad  \qquad \ \ \ \qquad g({\cal R}, z)=\frac {1}{1- Jz^Q}. \qquad  \ \ \qquad  \qquad  \qquad \qquad  \qquad \qed
$$
\renewcommand{\qedsymbol}{}
\end{pf}

\begin{cor}\label{cor: 2.6}
The topological entropy of $X(G_{\bigcirc}(N,J,Q))$ is equal to the negative logarithm of the positive root of
$$
z^Q+ \frac{N+1}{J}z - \frac{1}{J} =0.
$$

\end{cor}
\begin{pf}
See Corollary 2.4.
\end{pf}

We note another special case.

\begin{cor}\label{cor: 2.7}
$$
\zeta_{X(G_{\bigcirc}(N,1,2))}(z) =
\frac{2(1-z^2 + \sqrt{(1-z^2)^2 - 4Nz^2})}{(1 - z^2 - 2Nz + \sqrt{(1-z^2)^2 - 4Nz^2})^2},
\qquad N > 1.
$$
\end{cor}

\begin{cor}\label{cor: 2.8}
$$
h(X(G_{\bigcirc}(N,1,2))) =  \log 2 - \log(\sqrt{(N+1)^2 +4} -N - 1),
 \qquad N > 1.
$$
\end{cor}

As another example take for $G_{\circ}$  the Fischer automaton of the even system \cite [Section 6.1] {Ki}  
There are two vertices $v_{even}$ and $v_{odd}$, the labeling alphabet is $\{ 0,1 \}$,  and there are 
edges from  $v_{even}$  to $v_{odd}$ and from $v_{odd}$  to $v_{even}$  that are assigned by a
labeling map $\lambda_{even}$  the label $0$, and there is a loop at   $v_{even}$  that is assigned by the
labeling map $\lambda_{even}$  the label  $1$. $X (G(\lambda_{even}, v_{even}, (1,1)_{\gamma \in \Gamma}))$ is the Schr\"oder shift $Sch_{N}, N > 1$.

\begin{cor}\label{cor: 2.9}
$$
\zeta_{Sch_{N} }(z)=\frac{2(1+z)(1-z-z^2 + \sqrt{(1-z-z^2)^2 - 4Nz^2})}
{(1 -(2N+1)z- z^2  + \sqrt{(1-z-z^2)^2 - 4Nz^2})^2},
\qquad N > 1.
$$
\end{cor}
\begin{pf}
In this case by (2.1)
$$\ \ \
\qquad \qquad \qquad \qquad \qquad \qquad  g( {\cal R} , z)= \frac{1}{1 - z-z^2}. \qquad  \qquad\qquad  \qquad \qquad \qquad \qed
$$
\renewcommand{\qedsymbol}{}
\end{pf}

\begin{cor}\label{cor: 2.10}
$$
h(Sch_{N}) = \log 2  - \log (\sqrt{(N+2)^2 +4}-N-2), \qquad N > 1.
$$
\end{cor}
\begin{pf}
See Corollary 2.4.
\end{pf}

\begin{cor}\label{cor: 2.11}
\begin{multline*}
\zeta_{X(G(\lambda_{even}, v_{odd}, (1,1)_{\gamma \in \Gamma}))}(z) = \\
\frac{2(1+z)(1-z-z^2 + \sqrt{(1-z-z^2)^2 - 4Nz^2(1-z)^2})}
{(1 - (2N+1)z + (2N - 1)z^2 + \sqrt{(1-z-z^2)^2 - 4Nz^2(1-z)^2})^2},
\qquad N > 1.
\end{multline*}
\end{cor}
\begin{pf}
In this case by (2.1)
$$\ \ \
\qquad \qquad \qquad \qquad \qquad \qquad g({ \cal R} , z)= \frac{1-z}{1 - z-z^2}. \qquad  \qquad\qquad  \qquad \qquad \qquad \qed
$$
\renewcommand{\qedsymbol}{}
\end{pf}

\begin{cor}\label{cor: 2.12}
\begin{multline*}
h(X(G(\lambda_{even}, v_{odd}, (1,1)_{\gamma \in \Gamma}))) = \log 2 + \log N - \log (N+ 2 - \sqrt{(N+2)^2 -4N}), \\ N > 1,
\end{multline*}
\end{cor}
\begin{pf}
See Corollary 2.4.
\end{pf}

\section{Assigning to the elements of a finite set subsets of the set}

Let $\Psi$ be a map that assigns to a $\gamma \in \Gamma  $ a non-empty subset  $\Psi( \gamma) $ of $\Gamma $. We say that a permutation $ \pi$ of $\Gamma $  is a symmetry of $\Psi$ if
$$
\pi( \Psi( \gamma) ) = \Psi  ( \pi (\gamma)),   \qquad    \gamma \in \Gamma .
$$
We introduce an equivalence relation $ \sim$ into the set $\Gamma $ where for
 $\alpha, \beta \in \Gamma, \alpha\sim \beta$ means that for some $L \in \Bbb N$ there are $\gamma_{l}\in \Gamma, 0 \leq l \leq L, \gamma_{0} = \alpha, \gamma_{L} = \beta,$ such that  $\Psi( \gamma_{l-1}) = 
 \Psi( \gamma_{l})$ or there exists a symmetry $  \pi$ of $\Psi $ such that 
 $\pi( \gamma_{l-1}) =  \gamma_{l}, 0 \negthinspace < l \leq\negthinspace L$. 
 The set
 $$
 \Delta_{\Gamma} = \Psi^{-1}(\Gamma )
 $$
is an  $ \sim$-equivalence class. We denote by  $\Delta_{\setminus} $ the set of elements of $\Gamma $ that are    $ \sim$-equivalent to a $ \gamma \in \Gamma  $  such that
$$
 \Psi( \gamma ) = \Gamma \setminus \{ \gamma  \},
$$
and we denote by  $\Delta_{\bullet} $ the set of elements of $\Gamma $ that are   $ \sim$-equivalent to a $ \gamma \in \Gamma  $  such that
$$
 \Psi( \gamma ) =  \{ \gamma  \}.
$$
We set
$$
\Delta_{\circ} = \Gamma \setminus (   \Delta_{\Gamma}   \cup   \Delta_{\setminus}    \cup  \Delta_{\bullet} ),
$$
and we set
$$
N_{\Gamma} = \card ( \Delta_{\Gamma}),N_{\setminus} = \card ( \Delta_{\setminus}),N_{\bullet} = \card ( \Delta_{\bullet}),\
$$
$$
N_{\circ} = N - ( N_{\Gamma}   + N_{\setminus}  + N_{\bullet}   ).
$$
We denote the set of  $\sim$-equivalence classes in $\Delta_{\circ}$ by ${\cal A}$.

\begin{lem}\label{lem:3.1}
\begin{align*}
&\card ( \Psi( \alpha ) \cap  \Delta_{\Gamma}   ) =
 \card ( \Psi( \alpha^{\prime} ) \cap  \Delta_{\Gamma}   ) , \\
&\card ( \Psi( \alpha ) \cap  \Delta_{\setminus}  )= \card (\Psi(  \alpha^{\prime} ) \cap 
\Delta_{\setminus}  ),\\
&\card ( \Psi( \alpha ) \cap  \Delta_{\bullet}  )= \card ( \Psi( \alpha^{\prime} ) \cap  \Delta_{\bullet}  ), \qquad  \alpha,  \alpha^{\prime} \in A \in {\cal A}.
\end{align*}
\end{lem}
\begin{pf}
For a symmetry $\pi$ of $\Psi$,
\begin{align*}
\qquad \qquad \qquad&\card ( \Psi(\gamma ) \cap  \Delta_{\Gamma}   ) =
 \card ( \Psi(\pi(\gamma)) \cap  \Delta_{\Gamma}   ) , \\
&\card ( \Psi( \gamma ) \cap  \Delta_{\setminus}  )= \card (\Psi( \pi(\gamma ) )\cap 
\Delta_{\setminus}  ),\\
&\card ( \Psi( \gamma) \cap  \Delta_{\bullet}  )= \card ( \Psi(\pi( \gamma)) \cap  \Delta_{\bullet}  ), \qquad  \gamma \in \Gamma. \qquad  \qquad \qquad\qed
\end{align*}
\renewcommand{\qedsymbol}{}
\end{pf}

In view of  Lemma 3.1 we can introduce the notation
\begin{align*}
&K_{A}(\Gamma) =\card ( \Psi( \alpha ) \cap  \Delta_{\Gamma}   ),\\
&K_{A}(\setminus)=\card ( \Psi( \alpha ) \cap  \Delta_{\setminus}  ),\\
&K_{A}(\bullet)=\card ( \Psi(  \alpha) \cap  \Delta_{\bullet}  ), \qquad \alpha \in A \in {\cal A}.
\end{align*}

\begin{lem}\label{lem:3.2}
$$
\card ( \Psi( \alpha ) \cap  B  ) =\card ( \Psi( \alpha ^{\prime}) \cap B  ), \qquad
\alpha,  \alpha^{\prime} \in A \in {\cal A}, B \in {\cal A}.
$$
\end{lem}
\begin{pf}
For a symmetry $\pi$ of $\Psi$,
$$
 \quad \qquad\qquad\card ( \Psi( \gamma ) \cap  B  ) =\card ( \Psi(\pi(\gamma) \cap B  ), \qquad
\gamma \in \Gamma, B \in{\cal A}.  \qquad \qquad \quad \qed
$$ 
\renewcommand{\qedsymbol}{}
\end{pf}
In view of  Lemma 3.2 we can introduce the notation
$$
K_{A}(B)=\card ( \Psi(  \alpha) \cap  B  ), \qquad \alpha \in A \in {\cal A}, B \in{\cal A}.
$$

\section{Excluding a literal-uniform  set of words of length two from Dyck shifts}

We continue to consider  a mapping $\Psi$ that assigns to $ \gamma  \in  \Gamma $  a non-empty subset $\Psi(\gamma)$ of $\Gamma$.   Denote by 
 $X_{\Psi}$ the subshift that is obtained by removing from  the Dyck shift $D_{N}$ the words 
 $$\beta(+)\alpha(+), \qquad \alpha \in  \Gamma  ,\beta \in  \Gamma,  \beta \notin \Psi( \alpha ) .
 $$ 
Denote by ${\cal D}_{\gamma} $ the language of  words in the Dyck code that begin with $\gamma(-)$ and that are admissible for  $X_{\Psi}$. We set 

$$
\xi = 1 - \sum_{\gamma \in \Gamma} g( { \cal D}_{\gamma}).
$$

\begin{lem}\label{lem:4.1}
$$
g( { \cal D}_{\alpha}  , z) = z^2(1 + \frac{1}{\xi(z)} \sum_{ \beta \in \Psi( \alpha )}g( { \cal D}_{\beta}  , z)), \qquad \alpha \in \Gamma.
$$
\end{lem}
\begin{pf}
Apply the Sch\"utzenberger  method.
\end{pf}

\begin{lem}\label{lem:4.2}
Let $1 \leq K \leq N$, and let
$$
\card(\Psi(\gamma)   )= K, \qquad \gamma \in \Gamma.
$$
Then there exists for $\alpha,\beta \in \Gamma $ a length preserving bijection
$$
\eta(  \alpha,\beta ):   {\cal D}_{\alpha}  \to   {\cal D}_{\beta} .
$$
\end{lem}
\begin{pf}
Denote
$$
h(\gamma(-)) = 1,  h(\gamma(+) )= -1 \qquad     \gamma \in \Gamma ,
$$
and for a Dyck code word $(d_{q})_{1\leq q \leq Q}$ set
$$
H(d) = \max _{1\leq q < Q} \sum_{1\leq r \leq q} h(  d_{q}).
$$
For $H\in \Bbb N$ denote by $   {\cal D}_{H}( \alpha)  $ the set of $d \in  {\cal D}_{\alpha} $ such that $H(d) = H$.
Choose for  $\alpha,\beta  \in \Gamma$ bijections
$$
\chi(  \alpha,\beta   ): \Psi( \alpha)  \to  \Psi( \beta) .
$$
One constructs inductively length preserving bijections
$$
\eta_{H}(\alpha,\beta   ) : {\cal D}_{H}( \alpha)  \to   {\cal D}_{H}( \beta), \qquad H \in \Bbb N,
$$
where one sets
$$
\eta_{1}(\alpha,\beta   )( \alpha(-)  \alpha(+) ) =  \beta(-)   \beta(+), \qquad\alpha,\beta \in \Gamma.
$$
Assume that for $ \alpha,\beta  \in \Gamma $ and for $H\in \Bbb N $, the bijections
$$
\eta_{ \widetilde{H}}(\alpha,\beta   ) :\ {\cal D}_{ \widetilde{H}}( \alpha)  \to  {\cal D}_{ \widetilde{H}}( \beta),
 \qquad 1 \leq \widetilde{H}\leq H ,
$$
have been constructed. Then let $\eta_{H+1}(\alpha,\beta   ) $ be the map that carries a word
$$
\alpha(-)(   d^ {(l)})_{1 \leq l \leq L}\alpha(+) \in {\cal D}_{H+1}( \alpha) ,
$$
where
$$
   d^ {(l)}\in  {\cal D}^0 ,  H(d^{(l)})
 \leq H, \qquad1 \leq l \leq L \in \Bbb N ,
$$
with  $\gamma \in \Psi( \alpha)$ given by
$
d^{(L)} \in   {\cal D}_{\gamma}, 
$
into the word
$$
\beta(-)(   d^ {(l)})_{1\leq l < L}\eta_{ H(d^{(L)} )}(\gamma, \chi(  \alpha,\beta   )(\gamma  ))(d^{(L)})
\beta(+). 
$$
Set
$$
\qquad\qquad \eta(\alpha,\beta   ) (d) =\eta_{H}(\alpha,\beta   )(d), \qquad d \in {\cal D}_{H}( \alpha) ,  L \in \Bbb N,
\quad\alpha,\beta \in \Gamma. \qquad \qquad\qed
$$
\renewcommand{\qedsymbol}{}
\end{pf}

\begin{lem}\label{lem:4.3}
Let $1 \leq K \leq N$, and let 
$$
\card  ( \Psi( \gamma))= K, \qquad \gamma \in \Gamma.
$$
Then 
$$
g(  { \cal D}_{\alpha},z) =  \frac {1}{2N}(1 + (N-K)z^2-\sqrt{(1 + (N-K)z^2)^2-4Nz^2}\ ), \qquad\alpha \in \Gamma.
$$
\end{lem}
\begin{pf}
By Lemma 4.1 and Lemma 4.2
$$
 \quad \qquad \qquad \  \qquad g(   {\cal D}_{\alpha},z) = z^2 (1 + \frac{Kg(   {\cal D}_{\alpha},z) }{1 - N g(   {\cal D}_{\alpha},z) }), \qquad \alpha \in \Gamma. \qquad \quad \qquad \qquad\qed
$$
\renewcommand{\qedsymbol}{}
\end{pf}

\begin{lem}\label{lem:4.4}
$$ 
g(    {\cal D}_{\gamma},z) = \frac{z^2}{\xi(z)}, \qquad \gamma \in \Delta_{\Gamma}.
$$
\end{lem}
\begin{pf}
Apply the Sch\"utzenberger method.
\end{pf}

\begin{lem}\label{lem:4.5}
$$ 
g(  {\cal D}_{\gamma},z) = \frac{z^2}{\xi(z)+z^2}, \qquad \gamma \in \Delta_{\setminus}.
$$
\end{lem}
\begin{pf}
By Lemma 4.1
\begin{align*}   
   \qquad \qquad  \qquad g(   {\cal D}_{\gamma},z) &=   z^2(1 + \frac{1}{\xi(z)} \sum_{ \beta  \neq \gamma}g(   {\cal D}_{\beta}  , z))   \\
& =    z^2(1 + \frac{1}{\xi(z)} (1 - \xi (z)  -g(   {\cal D}_{\gamma},z)   )), \qquad  \gamma \in \Delta_{\setminus}.  \qquad  \ \qquad\qed
\end{align*}
\renewcommand{\qedsymbol}{}
\end{pf}

\begin{lem}\label{lem:4.6}
$$ 
g(    {\cal D}_{\gamma},z) = \frac{z^2\xi(z)}{\xi(z)-z^2}, \qquad \gamma \in \Delta_{\bullet}.
$$
\end{lem}
\begin{pf}
By Lemma 4.1
$$
\qquad \qquad\qquad \qquad g(   {\cal D}_{\gamma},z) = z^2(1+  \frac{g(   {\cal D}_{\gamma},z)}{\xi(z)}  ), \qquad  \gamma \in \Delta_{\bullet}.\qquad \qquad\qquad \qquad\qed
$$
\renewcommand{\qedsymbol}{}
\end{pf}

\begin{lem}\label{lem:4.7}
 $$
 g(   {\cal D}_{\alpha}) =  g(   {\cal D}_{\beta} ), \qquad \alpha,\beta \in A \in  {\cal A} .
 $$
\end{lem}
\begin{pf}
Apply Lemma 4.1. Also observe that a symmetry $\pi$ of $\Psi$ provides a length preserving bijection of ${\cal D} _{\gamma}$ onto ${\cal D} _{\pi(\gamma)},  \gamma \in \Gamma.$
\end{pf}

\begin{lem}\label{lem:4.8}
Let all elements of $\Delta_{\circ}$ be $\sim$-equivalent. Then
\begin{multline*}
g(    {\cal D}_{\gamma}, z ) =\\
\frac{z^2}{\xi(z) - K_{\Delta_{\circ}}( \Delta_{\circ} )z^2}\left( \xi (z)+   
\frac{K_{\Delta_{\circ}}( \Gamma) z^2}{\xi(z)} +  
\frac{K_{\Delta_{\circ}}( \bullet) z^2\xi(z)}{\xi(z)-z^2}  +  
\frac{ K_{\Delta_{\circ}}( \setminus )z^2}{\xi(z)+z^2}  \right)  , \\
\gamma\in   \Delta_{\circ}.
\end{multline*}
\end{lem}
\begin{pf}
In view of Lemma 4.7 one can introduce the notation
$$
g_B =g(    {\cal D}_{\beta} ) , \qquad \beta \in B \in  {\cal A},
$$
and rewrite then the formula of Lemma 4.1 as
\begin{align*}
&g_{A}(z)=  \\
&z^2\negthinspace\left(\negthinspace 1 + \frac{1}{\xi(z)}   \left (\negthinspace 
\frac{ K_{A}( \Gamma)  z^2}{\xi(z)} +  
 \frac{K_{A}(\bullet )z^2\xi(z)}{\xi(z)-z^2}  +  
\frac{K_{A}(\setminus )z^2}{\xi(z)+z^2} + \sum_{B \in  {\cal A}} K_{A}(B)g_B(z)\negthinspace \right) \negthinspace \right)\negthinspace, A\in {\cal A}.
\end{align*}
\end{pf}
 
\negthinspace

In case that $\Delta_{\circ}$ is empty or in case that all elements of $\Delta_{\circ} $ are $\sim$-equivalent one obtains from Lemma 4.4, Lemma 4.5, Lemma 4.6 and Lemma 4.8 an equation for $\xi$ that is at most of the fifth degree. Once solved, this equation yields then the generating functions, and conequently  the zeta function, from which, given $N$, the  topological enrtropy can be (numerically) determined. Unless Lemma 4.2 applies this equation will be of at least the third degree and for $N \leq 3$ this equation is of  at most the fourth degee. We determine here the zeta function  for the case of Lemma 4.2.

\begin{thm}\label{thm:4.9}
Let $1 \leq K \leq N$, and let
$$
\card(\Psi(\gamma)   )= K, \qquad \gamma \in \Gamma.
$$
Let a transition matrix $A_{\Psi}$ be given by
$$
A_{\Psi}(\beta, \alpha) = \begin{cases} 1, &\text {if $\beta \in \Psi (\alpha),$}\\
0, & \text {if $\beta \notin \Psi (\alpha), \qquad \alpha, \beta \in\Gamma.$}
 \end{cases}
$$
With
$$
g(z) =  \tfrac {1}{2}(1 + (N-K)z^2-\sqrt{(1 + (N-K)z^2)^2-4Nz^2}\ ),
$$
one has
$$
\zeta_{X_{\Psi}}
(z)= 
\frac{(1-Kz)(1-g(z))}{{\det ({\bold 1} - A_{\Psi}}z)(1-Nz-g(z))(1-Kz-g(z))}.
$$
\end{thm}
\begin{pf}
Set
$$
{\cal C}^{0} = \bigcup_{\alpha \in \Gamma}{ \cal D}_{\alpha} .
$$
and let ${\cal C}^{-}$ be  the code that contains the words that are concatenations of a word in 
${\cal C}^{0}$  and a (possibly empty) word in the symbols $\gamma(-), \gamma \in \Gamma$, and let 
$\cal C^{+}$ be  the code that contains the words in ${\cal L} (X_{\Psi})$ that are concatenations of a word 
in ${\cal C}^{0}$  and a (possibly empty) word in the symbols $\gamma(+), \gamma \in \Gamma $. 
Let $Y^-$ be the full shift with alphabet $\{\gamma(-): \gamma \in \Gamma\}$ and let $Y^+$ be the topological Markov shift with alphabet $\{\gamma(+): \gamma \in \Gamma\}$ and transition matrix 
$A_{\Psi}$. 
$X_{\Psi}$ and $Y^-$,$Y^+, {\cal C}^{-}$,${\cal C}^{0},{\cal C}^{+}$ satisfy  the relations (1.3 - 6).
It is 
$$
g({\cal C} ^-)(z) = \frac{g(z)}{1-Nz}.
$$
Denote by ${\cal L}_\alpha$ the language of  words that are admissible for $Y^+$ and that  have $\alpha\in \Gamma$ as their first symbol.  By Lemma 4.3 one has
\begin{align*}
g({\cal C}^+, z) = &g(z) + \sum_{\{\alpha, \beta \in \Gamma: \beta \in \Psi(\alpha) \}} g({\cal D}_{\beta}, z) 
g({\cal L}_{\alpha}, z) =\\
&g(z)\left(1 +  \tfrac{1}{N}\sum_{\{\alpha, \beta \in \Gamma: \beta \in \Psi(\alpha) \}} g({\cal L}_{\alpha}, z)\right) =\\
&g(z)\left(1 +  \tfrac{1}{N}\sum_{\alpha, \beta \in \Gamma, q \in \Bbb N }(A^{q}_{\Psi})_{\beta, \alpha}z^q\right)=\\
&g(z)\left(1+   \sum_{\ q \in \Bbb N }K^{q}z^q\right)=    
 \frac{g(z)}{1-Kz}.
\end{align*}
The theorem follows now from (1.7).
\end{pf}

\begin{cor}\label{cor: 4.10}
\begin{multline*}
h(X_{\Psi} ) =  \log N+\log(N-K)+ \log 2 \\
 -\log \left( \sqrt{(N^2 +K)^2 + 4(N-K)N^2} -N^2 -K \right).
\end{multline*}
\end{cor}
\begin{pf}
One checks that
$$
\limsup  \frac{1}{n}\log \card  \  P_n(X_{\Psi} ) =  \lim_{n \to \infty} \frac{1}{n}\log \card  \ 
{\cal  L}_n(X_{\Psi} ),
$$
and the corollary follows from Theorem 4.9.
\end{pf}

For excluding literal-uniform words of length two from the Motzkin shifts one obtains results of the same scope by the same methods.

\section{Excluding a literal-uniform set of words from \\ Dyck shifts and Motzkin shifts:  Two examples}
\begin{thm}\label{thm:5.1}
Let $X_{N}$ be the subshift that is obtained by excluding from the Dyck shift $D_{N}$ the words
$$
\gamma^{\prime}(+)  \gamma (+)
\gamma^{\prime \prime}(+), \quad\gamma, \gamma^{\prime} , 
\gamma^{\prime \prime}  \in \Gamma,    \gamma^{\prime}  \neq 
\gamma^{\prime \prime}.
$$
With 
$$
\rho(z) =\sqrt{(1-(N+1)z^2)^2  -4N^2z^4  }  ,
$$
one has
\begin{multline*}
 \zeta_{X_{N}}(z)= \\
\dfrac{ 2 (  1- (N-1)z^2 + \rho(z)   )  }      
    { ( 1 - 2Nz -(N-1)z^2     +\rho(z)  ) (     1 -2z +(N+1)z^2 - 2N^2z^3+ \rho(z)   )      } \times\\
\frac {1}{(1-z)^{N-1}(1-z^2)^{\frac{N(N-1)}{2}} }.
\end{multline*}
\end{thm}
\begin{pf}
Take as the code ${\cal C}^-$ the circular  code that contains the words in ${\cal L}(X_{N})$ that are concatenations of  a word in the symbols $\gamma(-), \gamma \in \Gamma$ and of a word in the Dyck code. Take as the code ${\cal C}^0$ the circular code that contains the words in the Dyck code that are admissible for $X_{N}$. It is
$$
g({\cal C}^-,z) = \frac{g({\cal C}^0,z)}{1-Nz}.
$$
Take as the code  ${\cal C}^+$ the circular code that contains the words in ${\cal L}(X_{N})$ that are conactenations of a word in the Dyck code  and of a word in the symbols  $\gamma(+), \gamma \in \Gamma$.  Consider the code
$$
{\cal K}= {\cal C}^0 \setminus \{ \gamma(-) \gamma(+) : \gamma \in \Gamma\}.
$$
It is
$$
g({\cal C}^+,z) = \frac {N^2z^3 + g({\cal K} , z)}{1-z}.
$$
The Sch\"utzenberger method gives
$$
g({\cal K}, z) = \frac{N^2z^4 + z^2g({\cal K}, z) }{1-g({\cal K}, z) - Nz^2},
$$
leading to
$$
g({\cal K}, z) = \tfrac{1}{2}(1 - (N+1)z^2-\rho (z)),
$$
and
\begin{align*}
&g({\cal C}^- , z) =  \frac{1 + (N-1)z^2-\rho (z)}{2(1-Nz)}, \tag {5.1}\\
&g(({\cal C}^0 , z) =  \tfrac{1}{2}(1 + (N-1)z^2-\rho (z)), \tag{ 5.2}\\
&g({\cal C}^+,z) = \frac {2N^2z^3 +1 - (N+1)z^2-\rho (z)}{2(1-z)}. \tag {5.3}
\end{align*}
Take as $Y^-$ the full shift with  symbols  $\gamma(-), \gamma \in \Gamma$, and take as $Y^+$ the set of points of period two in the full shift with alphabet $\{\gamma(+): \gamma \in \Gamma\}$.
$X_{N}, {\cal C}^-,{\cal C}^0 ,{\cal C}^+,Y^-,Y^+$ satisfy the relations (1.3 - 5) and the theorem follows from (1.6) together with with (5.1 - 3).
\end{pf}

\begin{prop}\label{prop:5.2}
Let $X^{(N,1)}$ be the subshift that is obtained by excluding from the Motzkin shift $ M_{N}$ the words
$$
\gamma (+) \gamma^{\prime} (-) ,\gamma (+) \gamma^{\prime} (+) , \quad  \gamma, \gamma^{\prime}\in \Gamma,
$$
as well as the words
$$
\bold 1 \bold 1 \gamma (+), 
\gamma(-) \bold 1 \gamma(+), \quad \gamma  \in \Gamma.
$$
With
$$
\rho(z) = \tfrac{1}{2}(1 - z - \sqrt{(1 - z)(1 - z - 4Nz^3)}),
$$
one has
$$
\zeta  _{X^{(N,1)}}(z) =
\frac{1-z -g(z)}{(1-(N+1)z - g(z)) ((1 - Nz^2)(1-z)-g(z))}.
$$
\end{prop}
\begin{pf}
Denote  by ${\cal C}$ the circular code that contains the words that are obtained by putting the symbol 
$\bold 1$ at the end of a word in the Motzkin code that is admissible for  $X^{(N,1)}$. Take as the code ${\cal C}^{-}$ the circular code that contains the words  that are concatenations of  a word in the symbols of $ \{    \gamma (-): \gamma \in \Gamma\}  \cup \{ {\bold 1}\}$
and a word in ${\cal C}$. It is
\begin{align*}
g({\cal C}^- , z) =  \frac{g({\cal C} , z)}{1-(N+1)z}. \tag {5.4}
\end{align*}
Take as the code ${\cal C}^{0}$ the circular code that contains the words that are obtained by letting a  strings of $\bold 1$'s follow a word in  ${\cal C}$. It is
\begin{align*}
g({\cal C}^0 , z) =  \frac{g({\cal C} , z)}{1-z}, \tag {5.5}
\end{align*}
Take as the code
${\cal C}^{+}$ the circular code that contains the words in ${\cal L}(X^{(N,1)})$ that are concatenations of a word in ${\cal C}$ and a word in the symbols of $ \{\gamma(+): \gamma \in \Gamma \} \cup \{\bold 1\}$.  It is
\begin{align*}
g({\cal C}^+,z) =  \frac{g({\cal C} , z)}{(1-Nz^2)(1-z)} \tag {5.6}.
\end{align*}
The Sch\"utzenberger method gives
$$
g({\cal C}, z) = \frac {Nz^3(1-z)}{1-z - g({\cal C}, z)},
$$
and it follows that 
$$
g({\cal C}, z) = \tfrac{1}{2}(1 - (N+1)z^2 -\rho(z)).
$$
Take as $Y^-$ the full shift with alphabet $\{\gamma(-): \gamma \in \Gamma\}\cup \{{\bold 1}\}$, and take as $Y^+$ the subshift of finite type that is obtained  by removing the words $\bold 1\bold 1   \gamma(+), \gamma \in \Gamma $ and the words $\gamma (+) \gamma^{\prime} (+) , \quad  \gamma, \gamma^{\prime}\in \Gamma$, from the full shift with alphabet $\{\gamma(+): \gamma \in \Gamma\}\cup \{{\bold 1}\}$.  $X^{(N,1)}, {\cal C}^-,{\cal C}^0,{\cal C}^+, Y^-,Y^+$ satisfy the relations (1.4 - 6)  and  the theorem follows from (1.7).
\end{pf}

\section{Excluding a literal-uniform set of words \\ from $ X(G_{\bigcirc}, N, J, 1)$: An example }

We return to the graph $G_{\bigcirc}(N, J, 1), J \in \Bbb N,$ of section 2. We denote the set of edges of  $G_{\bigcirc}( N,J, 1)$ by $\Omega$.

\begin{thm}\label{thm:6.1}
Let $J , K\in{ \Bbb N},1 \leq K, L \leq J $. Let $ \Xi_\Omega$ be a mapping that  assigns to $\omega \in \Omega$ a non-empty subset $\Xi (\omega)$ of $\Omega$, such that
\begin{align*}
\card (\Xi_\Omega (\omega)) = K, \quad \omega \in \Omega.  \tag {6.1}
\end{align*}
and let
$ \Xi_\Gamma$ be a mapping that  assigns to $\gamma \in \Gamma$ a non-empty subset 
$\Xi (\gamma)$ of $\Omega$, such that
\begin{align*}
\card (\Xi_\Gamma (\gamma)) = L, \quad \gamma \in \Gamma.  \tag {6.2}
\end{align*}
Let $X_{N, \Xi_\Omega, \Xi_\Gamma} $ be  \thinspace the \thinspace subshift   \thinspace that  \thinspace is  \thinspace  obtained  \thinspace by \thinspace excluding  \thinspace from  \thinspace  the subshift $ X(G_{\bigcirc}( N, J, 1))$ the words
$$
\gamma (+) \gamma^{\prime} (-) ,\gamma (+) \gamma^{\prime} (+) , \quad  \gamma, \gamma^{\prime}\in \Gamma,
$$
as well as the words
$$
\omega \omega^{\prime}\gamma (+), \quad \omega , \omega^{\prime}  \in \Omega, \gamma \in \Gamma,
$$
$$
\gamma (-) \omega\gamma(+), \quad\gamma\in \Gamma, \omega \in \Omega,
$$
and the words
$$
\omega \gamma (+)\omega^{\prime}, \quad  \gamma \in \Gamma, \omega , \omega^{\prime}  \in \Omega,  \omega^{\prime}  
\notin \Xi_\Omega (\omega), 
$$
$$
\gamma (-) \gamma (+)  \omega, \quad  \gamma \in \Gamma,  \omega \in \Omega,   \omega \notin   \Xi_\Gamma(\gamma).
$$
Assume that the  transition matrix $A$ that is given given by
$$
A(\omega, \omega^{\prime})= \begin{cases} 1, &\text {for $ \omega^{\prime} \in \Xi_\Omega (\omega),$}\\
0, & \text {for $ \omega^{\prime} \notin \Xi_\Omega(\omega), \quad \omega , \omega^{\prime}\in \Omega,$}
 \end{cases}
$$
is irreducible. With the period $\pi$ of the matrix $A$ define  a polynomial $q$ by
$$
\det ({ \bold 1} -  Az  ) =  q(z)(1 - K^{\pi}z^{\pi }).
$$
With
$$
g(z) = \tfrac{1}{2}(1\negthinspace - \negthinspace Jz \negthinspace+ \negthinspace
 N(L \negthinspace-K \negthinspace)z^3 \negthinspace- \negthinspace\sqrt{(1\negthinspace \negthinspace- \negthinspace Jz\negthinspace + \negthinspace N(L\negthinspace-\negthinspace K)z^3)^2\negthinspace -\negthinspace 4NLz^3(1\negthinspace -\negthinspace Jz)}),
$$
one has
\begin{multline*}
\zeta  _{X_{N, \Xi_\Omega, \Xi_\Gamma} }(z) =
\\
\frac{1-Jz -g(z)}{(1-(J+N)z - g(z)) ((1 - NKz^2)(1-Jz)-g(z))q(Nz^2)\sum_{0\leq r< \pi}N^rK^rz^{2r}}.
\end{multline*}
\end{thm}
\begin{pf}
For $\omega \in \Omega$ denote  by ${\cal C}_{\omega}$ the circular code that contains the words that are obtained by putting the symbol  $\omega$ at the end of a word in ${\cal L} (X_{N, \Xi_\Omega, \Xi_\Gamma} )$ whose 
$\lambda_{\bigcirc}$-label sequence is in the Motzkin code. For $\omega \in \Omega$ denote by 
${\cal C}^+_{\omega}$  the code that contains the words in ${\cal L} (X_{N, \Xi_\Omega, \Xi_\Gamma} )$ that are concatenations of a word in  ${\cal C}_{\omega}$ and of a word that takes its symbols from $\Omega \cup \{\gamma(+): \gamma \in \Gamma \}$, and that has its last symbol in $ \Omega$. By (6.1)
\begin{align*}
 g({\cal C}^{+}_{\omega} , z)  =
\frac {g( {\cal C}_{\omega} , z)}{(1-NKz^2)(1-Jz)}.  \tag {6.3}
\end{align*}
Set
$$
{\cal C}= \bigcup_{\omega \in \Omega }{\cal C}_{\omega}. 
$$
As the code ${\cal C}^{-}$ take the the circular code that contains the words  that are concatenations of a word that takes its  symbols from $\Omega \cup \{\gamma(-): \gamma \in \Gamma \}$ and of a word in  ${\cal C}$. It is
\begin{align*}
g( {\cal C}^-,z) =
\frac{g({\cal C},z)}  { 1-(N+J)z }.\tag {6.4}
\end{align*}
As the code ${\cal C}^{0}$ take the circular code that contains the words that are concatenations of a word in  ${\cal C}$ and of a word that takes its  symbols from $\Omega$. It is
\begin{align*}
g( {\cal C}^0,z) =
\frac{g({\cal C},z)}  { 1-Jz }. \tag {6.5}
\end{align*}
As the code ${\cal C}^{+}$ take the circular code that conatains the words in ${\cal L} (X_{N, \Xi_\Omega, \Xi_\Gamma} )$ that are concatenations of a word in  ${\cal C}$ and of a word that takes its symbols from
$\Omega \cup \{\gamma(+): \gamma \in \Gamma \}$, and that has its last symbol in $ \Omega$. It is
$$
{\cal C}^{+} = \bigcup_{\omega \in \Omega }{\cal C}^+_{\omega},
$$ 
and it follows from (6.3) that 
\begin{align*}
g({\cal C}^{+} , z) =  
\sum_{\omega \in \Omega } 
\frac {g({\cal C}_{\omega} , z)}{(1-NKz^2)(1-Jz)} = \frac{g({\cal C} , z) }{(1-NKz^2)(1-Jz)}. \tag {6.6}
\end{align*}
As $Y^-$ take the full shift with alphabet  $\Omega \cup \{\gamma(-): \gamma \in \Gamma \}$.  It is
\begin{align*}
\zeta _{Y^-}(z)= \frac{1}{1 - (N+J)z}. \tag {6.7}
\end{align*}
As $Y^+$ take the subshift of finite type that is obtained by excluding from the full shift with alphabet  
$\Omega \cup \{\gamma(+): \gamma \in \Gamma \}$ the words
$$
\gamma (+) \gamma^{\prime} (+) , \quad  \gamma, \gamma^{\prime}\in \Gamma,
$$
as well as the words
$$
\omega \omega^{\prime}\gamma (+), \quad \omega , \omega^{\prime}  \in \Omega, \gamma \in \Gamma,
$$
and the words
$$
\omega \gamma (+)\omega^{\prime}, \quad \omega , \omega^{\prime}  \in \Omega,  \omega^{\prime}  
\notin  \Xi_\Omega (\omega),  \gamma \in \Gamma.
$$
It is
\begin{align*}
\zeta_{Y^+}(z) = \frac {1}{
(1-Jz)
 q(Nz^2)
 (1-N^\pi K^\pi z^{2\pi})} .   \tag {6.8}
\end{align*}

Further choose a $\tilde {\gamma} \in \Gamma$ and
denote by $\tilde{ {\cal C}}^+_{\omega}$ the circular code that contains the words in  ${\cal L} (X_{N, \Xi_\Omega, \Xi_\Gamma} )$ that are concatenations of a word in ${\cal C}_{\omega}$ and of a word that takes its symbols from $\Omega \cup \{\tilde{\gamma} (+) \}$, and that has its last symbol in $ \Omega$. Again using (6.1) one has
\begin{align*}
 g(\tilde{{\cal C}}^{+}_{\omega} , z)  =
\frac {g(  {\cal C}_{\omega} , z)}{(1-Kz^2)(1-Jz)}.  \tag {6.9}
\end{align*}
Denote by $\tilde{{\cal C}}^+$ the circular code that contains the words in  ${\cal L} (X_{N, \Xi_\Omega, \Xi_\Gamma} )$ that are concatenations of a word in $ {\cal C}$ and of a word that takes its symbols from $\Omega \cup \{\tilde{\gamma} (+) \}$, and that has its last symbol in $ \Omega$. It is
$$
\tilde{ {\cal C}}^{+} = 
\bigcup_{\omega \in \Omega }
\tilde{{\cal C}}^+_{\omega},
$$ 
and by (6.8) then
\begin{align*}
g( \tilde{{\cal C}}^{+} , z)  = \frac{g( {\cal C}, z) }{(1-Kz^2)(1-Jz)}. \tag {6.10}
\end{align*}
Let $\tilde{Y}^+$ be the subshift of finite type that is obtained by excluding from the full shift with alphabet  
$\Omega  \cup \{ \tilde {\gamma}(+) \}$ the word
$
\tilde {\gamma}(+)\tilde {\gamma}(+)
$
as well as the words
$$
\omega \omega^{\prime}\tilde {\gamma}(+), \quad \omega , \omega^{\prime}  \in \Omega, $$
and the words
$$
\omega \tilde {\gamma}(+)\omega^{\prime}, \quad \omega , \omega^{\prime}  \in \Omega,  \omega^{\prime}  
\notin  \Xi_\Omega (\omega).
$$
It is
\begin{align*}
\zeta_{\tilde{Y}^+}(z) = \frac {1}{(1-Jz) q(z^2)(1- K^\pi z^{2\pi})}. \tag {6.11}
\end{align*}
Finally let 
 $\tilde{Y}$ be the subshift of finite type that is obtained by removing from the full shift with alphabet  
$\Omega \cup \{\gamma(-): \gamma \in \Gamma \} \cup \{ \tilde {\gamma}(+) \}$ the word
$
\tilde {\gamma}(+)\tilde {\gamma}(+)
$
and the words
$$
\tilde {\gamma}(+) \gamma (-)  , \quad  \gamma\in \Gamma,
$$
as well as the words
$$
\omega \omega^{\prime}\tilde {\gamma}(+), \quad \omega , \omega^{\prime}  \in \Omega, 
$$
$$
\gamma (-) \omega\tilde {\gamma}(+),, \quad\gamma ,\gamma^{\prime}\in \Gamma, \omega \in \Omega,
$$
and the words
$$
\omega \tilde {\gamma}(+)\omega^{\prime}, \quad \omega , \omega^{\prime}  \in \Omega,  \omega^{\prime}  
\notin \Xi_\Omega (\omega),
$$
$$
\gamma (-)  \tilde {\gamma}(+)\omega, \quad \gamma\in \Gamma,\omega \in \Omega,\omega
\notin \Xi_\Gamma (\gamma).
$$
It is by (6.1) and (6.2)
\begin{multline}
\zeta_{\tilde{Y}}(z) = \dfrac {1}{((1-Kz^2)(1-Jz) - Nz(1 + (L- K)z^2))q(z^2)\sum_{0\leq r< \pi}K^rz^{2r}}. \tag {6.12}
\end{multline}
There is a bijective correspondence between the set $P( \tilde{Y} )$ and the set
$$
P(Y^-) \cup P( Y( {\cal C}^- ) \cup Y( \tilde{ {\cal C}}^{+})) \cup P( \tilde{Y}^+ )
$$
 that preserves periods. To describe this correspondence we set
$$
\eta (\tilde {\gamma}(+)) = -1,
$$
and
$$
\eta (\gamma (-)) = 1, \quad \gamma \in \Gamma,\qquad \eta (\omega) = 0,\quad \omega \in \Omega, 
$$
and we denote for $p \in  P( \tilde{Y})$  by $ {\cal I} (p)$ the set of $i \in \Bbb Z$ such that $p_{i}=  \tilde {\gamma}(+) $, and such that the set
$$
 {\cal J}(p, i) = \{ j \in {\Bbb Z} : j < i, \sum_{j\leq k < i} \eta(p_{k})= 0 \}
$$
is not empty. We set
$$
j(p, i)= \max  {\cal J} (p, i), \quad p \in P( \tilde{Y} ), i \in {\cal I}(p).
$$
Here
$$p_{j(p, i)} \in \{\gamma(-) :  \gamma \in \Gamma \}.
$$
We set also
$$
\vartheta (\gamma(-) ) =  \gamma(+) , \quad \gamma \in \Gamma.
$$
The correspondence arises by simultaneously replacing in every point $p \in P( \tilde{Y} )$ for all  $ i \in  {\cal I}(p)$ the symbol $p_{i} = \tilde{\gamma}(+)$ by the symbol $\vartheta (p_{j(p, i) }) $. One has
$$
P(Y( {\cal C}^-)) \cap  P(Y( {\cal C}^+) ) = P(Y( {\cal C}^0) ),
$$
and the intersection of $Y^-$ and $Y^-$ is the full shift with alphabet $\Omega$.
By the existence of the correspondence we can therefore derive, using  (6.4 - 8) and (6.10 - 12), the equation
\begin{multline*}
\frac{1 - Jz - g ({\cal C}, z)}{(1 - (N+J)z-g( {\cal C}. z) )((1-Kz^2)(1-Jz)-g( {\cal C}. z)  )} = \\
\frac {1}{(1-Kz^2)(1-Jz) - Nz(1 + (L - K)z^2)},
\end{multline*}
which gives
\begin{multline}
g( {\cal C}, z) =  \\  \tfrac{1}{2}(1\negthinspace - \negthinspace Jz \negthinspace+ \negthinspace
 N(L \negthinspace- \negthinspace K)z^3 \negthinspace- \negthinspace\sqrt{(1\negthinspace \negthinspace- \negthinspace Jz\negthinspace + \negthinspace N(L \negthinspace-\negthinspace K)z^3)^2\negthinspace -\negthinspace 4NLz^3(1\negthinspace -\negthinspace Jz)}).
 \tag {6.13}
 \end{multline}
$X_{N, \Xi_\Omega, \Xi_\Gamma}  ,{\cal C}^-,{\cal C}^0, {\cal C}^+, Y^-,Y^+$ satisfy the relations (1.3 - 6) 
and  the  theorem follows from (1.7) and (6.13) by the use of (6.4 - 7).
\end{pf}

We note that Proposition 5.2 is the special case J = 1 of Theorem 6.1. Besides this case and the case that  $K = 1$, note also the case that $L=J$, and the case  that $J = N$ and $L = K$.

\end{document}